\documentclass[a4paper,10pt]{article}
\usepackage{amsmath,amssymb}
\usepackage[frenchb]{babel}
\usepackage[latin1]{inputenc}
\usepackage[dvips]{graphicx}

\newcommand{\x}{\mathsf{x}}
\newcommand{\BB}{\mathsf{B}}
\newcommand{\racine}{\#}
\newcommand{\extr}{\operatorname{Ext}(T)}

\newtheorem{theorem}{Théorème}[section] 
\newtheorem{proposition}[theorem]{Proposition}

\newtheorem{lemma}[theorem]{Lemme}

\newenvironment{proof}{\begin{trivlist}\item{\bf{Preuve.}}}
  {\hfill\rule{2mm}{2mm}\end{trivlist}}

\title{Une base symétrique de l'algèbre\\ des coinvariants quasi-symétriques}
\author{Frédéric Chapoton}
\date{\today}

\begin{document}

\maketitle

\begin{abstract}
  On décrit une nouvelle base de l'algèbre des coinvariants
  quasi-symétriques, qui est stable par l'involution naturelle et
  indexée par les triangulations d'un polygone régulier.
\end{abstract}

\setcounter{section}{-1}

\section{Introduction}

L'algèbre des coinvariants est un objet classique associé à chaque
groupe de Coxeter fini $W$ \cite{steinberg}. Cette algèbre est définie
comme le quotient de l'algèbre des polynômes sur l'espace vectoriel
sur lequel $W$ agit par réflexions, par l'idéal homogène engendré par
les polynômes invariants homogènes non constants. Le quotient est une
algèbre graduée de dimension finie donnée par l'ordre de $W$.

Dans le cas du groupe symétrique sur $n+1$ lettres, on peut expliciter
cette construction comme le quotient de l'algèbre des polynômes en
$\x_1,\dots,\x_{n+1}$ par l'idéal homogène engendré par les fonctions
symétriques élémentaires. Plus récemment, une notion plus faible que
celle de polynôme symétrique est apparue \cite{gessel}. Un polynôme
est dit quasi-symétrique si, pour toute suite d'exposants
$(m_1,\dots,m_k)$ fixée, tous les monômes $\x_{i_1}^{m_1}\dots
\x_{i_k}^{m_k}$ pour une suite \emph{croissante} d'indices
$i_1<i_2<\dots<i_k$ ont le même coefficient. En particulier, les
polynômes symétriques sont aussi quasi-symétriques.

Dans \cite{ba,bba}, la notion d'algèbre coinvariante quasi-symétrique
a été introduite et étudiée. Elle est définie comme le quotient de
l'algèbre des polynômes en $\x_1,\dots,\x_{n+1}$ par l'idéal homogène
engendré par les polynômes quasi-symétriques homogènes non constants.
C'est une algèbre graduée. Il est démontré dans \cite{bba} que cette
algèbre est de dimension finie, donnée par le nombre de Catalan
$c_{n+1}$. La preuve est la construction explicite d'un ensemble de
monômes indexés par les chemins de Dyck de longueur $2n+2$, dont les
images dans le quotient forment une base de l'algèbre des coinvariants
quasi-symétriques.

Dans le cas des coinvariants usuels, le groupe de Coxeter $W$ agit par
automorphismes sur le quotient et on obtient une décomposition
intéressante du module régulier. Dans le cas des coinvariants
quasi-symétriques, le seul automorphisme de la situation est le
renversement qui envoie $\x_i$ sur $\x_{n+2-i}$. Cette involution
préserve l'idéal des fonctions quasi-symétriques sans terme constant
et passe donc au quotient.

La motivation initiale de cet article est le fait que l'action de
cette involution semble difficile à décrire dans la base des monômes
associés aux chemins de Dyck. On construit donc une nouvelle base,
dans laquelle l'involution agit par permutation. Cette base est formée
de polynômes dont le terme dominant pour l'ordre naturel sur les
variables $\x_1,\dots,\x_{n+1}$ redonne les monômes associés aux
chemins de Dyck.

Il apparaît que l'ensemble naturel d'indexation de cette nouvelle base
est non pas l'ensemble des chemins de Dyck, mais celui des
triangulations d'un polygone régulier. Cet ensemble joue un rôle
primordial dans la théorie des algèbres à grappes de Fomin et
Zelevinsky \cite{cluster2}. Cet article donne donc un premier
rapprochement entre les algèbres à grappes et les fonctions
quasi-symétriques. Il se trouve que la construction de la base indexée
par les triangulations passe par le choix d'une triangulation de base
en forme d'éventail. Dans le cadre de la théorie des algèbres à
grappes, ce choix correspond au carquois équi-orienté de type $A_n$,
voir par exemple \cite{CCS}.

L'article est organisé comme suit. On commence par définir une
bijection \textit{ad hoc} entre triangulations et chemins de Dyck.
Ensuite on montre que, par cette bijection, le monôme dominant du
polynôme associé à une triangulation est le monôme associé au chemin
de Dyck correspondant, ce qui entraîne immédiatement le résultat
principal.
 
\section{Bijection}

\begin{figure}
  \begin{center}
    \includegraphics[width=8cm]{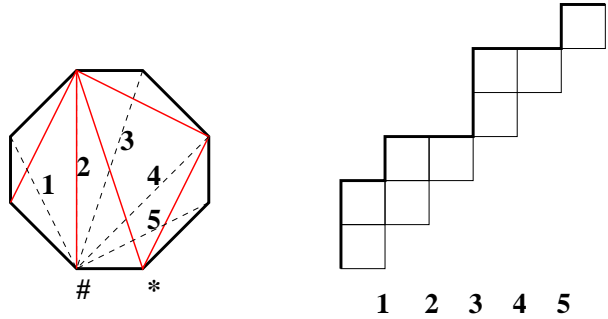}
  \end{center}
  \caption{Exemple pour la bijection}
  \label{fig:bijection}
\end{figure}

Soit $n$ un entier positif ou nul. On définit dans cette section une
bijection entre
\begin{enumerate}
\item les triangulations d'un polygone régulier à $n+3$ cotés,
\item les chemins de Dyck de longueur $2n+2$. 
\end{enumerate}
Il est bien connu que ces deux ensembles ont pour cardinal le nombre
de Catalan 
\begin{equation}
  c_{n+1}=\frac{1}{n+2}\binom{2n+2}{n+1}.
\end{equation}
Par définition, un chemin de Dyck est une suite de pas verticaux
(``montées'') et horizontaux (``descentes'') qui reste au dessus de la
diagonale, voir la partie droite de la figure \ref{fig:bijection}.

La bijection est illustrée par un exemple dans la figure
\ref{fig:bijection}.

Avant toute chose, on fixe une triangulation de base en forme
d'éventail, c'est-à-dire formée par toutes les diagonales contenant un
sommet choisi, noté $\racine$. On dessine cette triangulation avec le
sommet commun à toutes les diagonales placé en bas. Les diagonales de
cette triangulation de base seront dites ``négatives'' et numérotées
de $1$ à $n$ de gauche à droite. Les diagonales qui n'interviennent
pas dans la triangulation de base sont dites ``positives''. On
numérote aussi de $1$ à $n$ les sommets aux extrémités des diagonales
négatives.

\smallskip

On associe alors un chemin de Dyck $D(T)$ à chaque triangulation $T$,
par récurrence sur $n$. Pour $n=0$, à la seule triangulation du
polygone à trois cotés est associée le seul chemin de Dyck de longueur
$2$.

Si $n$ est non nul, on regarde le sommet $*$ du polygone placé à
droite du sommet $\racine$ dans le sens trigonométrique. On distingue
deux cas.

Si le sommet $*$ participe à un seul triangle de la triangulation $T$
\textit{i.e.} n'est contenu dans aucune diagonale de $T$, on lui
associe le chemin de Dyck obtenu en encadrant par une montée et une
descente le chemin de Dyck $D(T')$ associé à la triangulation $T'$ du
polygone à $n+2$ cotés qui est définie comme $T$ moins le triangle
adjacent à $*$. Le sommet distingué $\racine$ de $T'$ est celui de
$T$.

\begin{figure}
  \begin{center}
    \includegraphics[width=10cm]{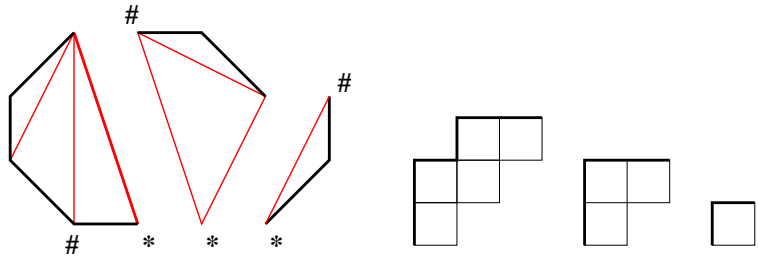}
  \end{center}
  \caption{Décomposition en morceaux}
  \label{fig:morceaux}
\end{figure}

Si le sommet $*$ participe à plusieurs triangles, on découpe la
triangulation en autant de morceaux (le long des diagonales contenant
$*$), voir la figure \ref{fig:morceaux}. Le sommet $*$ donne un sommet
dans chacun de ces morceaux. On prend dans chacun des morceaux le
sommet à gauche de $*$ comme sommet distingué $\racine$. Par
récurrence, on associe un chemin de Dyck à chacun des morceaux et on
les concatène dans l'ordre des morceaux induit par l'ordre
de gauche à droite au voisinage du sommet $*$ dans $T$, voir les
figures \ref{fig:bijection} et \ref{fig:morceaux}.

C'est clairement une bijection. La bijection inverse est aussi définie
par récurrence sur $n$. On décompose un chemin de Dyck réductible pour
la concaténation en ses composantes irréductibles et on recompose une
triangulation par juxtaposition. Pour les chemins de Dyck
irréductibles, on enlève une montée et une descente, on obtient une
triangulation par récurrence et on rajoute un triangle.

\begin{lemma}
  \label{premierspas}
  Le nombre de pas verticaux initiaux du chemin de Dyck $D(T)$ est le
  nombre de diagonales négatives dans $T$ plus $1$.
\end{lemma}

\begin{proof}
  La preuve se fait par récurrence. L'énoncé est vrai pour $n=0$. On
  distingue deux cas comme dans la définition de la bijection. Dans le
  cas où $*$ est dans une seule diagonale, les deux quantités
  augmentent de $1$. Dans l'autre cas, les deux quantités sont inchangées.
\end{proof}

\section{Polynômes}

On associe à chaque diagonale un polynôme en les variables
$\{\x_1,\dots,\x_{n+1}\}$ comme suit. On associe la constante $1$ aux
diagonales négatives. Chaque diagonale positive coupe un ensemble de
diagonales négatives consécutives de $i$ à $j$. En fait, ceci donne
une bijection entre les diagonales positives et les segments de
$\{1,\dots,n\}$. On peut donc parler de la diagonale positive $(i,j)$,
à qui on associe alors la somme des $\x_{k}-\x_{k+1}$ pour
$k=i,\dots,j$ soit $\x_{i}-\x_{j+1}$.

On associe alors à chaque triangulation $T$ le produit $\BB_T$ des
polynômes associés à ses diagonales. Dans l'exemple de la figure
\ref{fig:bijection}, on obtient
\begin{equation}
  \BB_T=(\x_1-\x_2)(1)(\x_3-\x_4)(\x_3-\x_6)(\x_5-\x_6).
\end{equation}

\smallskip

Par ailleurs, comme dans \cite{bba,ba}, on associe un monôme $M_D$ en
$\{\x_1,\dots,\x_n\}$ à chaque chemin de Dyck $D$. On représente un
chemin de Dyck par une suite de pas d'une unité vers le haut
(``montée'') ou vers la droite (``descente'') dans une grille. On
numérote les colonnes internes de la grille de $1$ à $n$, voir la
partie droite de la figure \ref{fig:bijection}. On convient que chaque
pas vertical d'indice $i$ correspond à la variable $\x_{i}$. Le monôme
$M_D$ est alors le produit des contributions des pas verticaux.  Dans
l'exemple de la figure \ref{fig:bijection}, on obtient
\begin{equation}
  M_{D(T)}=\x_1\x_3\x_3\x_5.  
\end{equation}

On définit un ordre sur les monômes en ordonnant les variables par
\begin{equation}
  \x_1 \ggg \x_2 \ggg \dots \ggg \x_{n+1}.
\end{equation}

Le monôme dominant d'un polynôme pour cet ordre est celui où
intervient la plus grande puissance de $\x_1$, puis en cas d'ambiguïté
la plus grand puissance de $\x_2$ et ainsi de suite.

\begin{proposition}
  \label{correspond}
  Le monôme dominant du polynôme $\BB_T$ associé à une triangulation
  $T$ est le monôme $M_{D(T)}$ associé au chemin de Dyck $D(T)$
  correspondant à $T$ via la bijection ci-dessus.
\end{proposition}

\begin{proof}
  Par récurrence sur $n$. La proposition est vraie pour $n=0$.
  Soit donc $n$ non nul. On distingue deux cas.

  Supposons d'abord que le sommet $*$ participe à un seul triangle de
  la triangulation $T$. Alors la triangulation $T$ contient la
  diagonale négative $n$. Le polynôme $\BB_T$ ne fait donc pas
  intervenir $\x_{n+1}$ et est égal au polynôme $\BB_{T'}$ associé à
  la triangulation raccourcie en $*$. De même, le chemin de Dyck
  $D(T)$ est obtenu par concaténation d'une montée, du chemin de Dyck
  $D(T')$ et d'une descente. Donc le monôme associé à $D(T)$ est le
  même que celui associé au chemin $D(T')$. On conclut par hypothèse
  de récurrence que le monôme dominant de $\BB_T$ est $M_{D(T)}$.

  Supposons maintenant que le sommet $*$ participe à plusieurs
  triangles de $T$. Soit $\extr$ l'ensemble des nombres $k$ dans
  $\{1,\dots,n\}$ tels que la diagonale négative $k$ partage un sommet
  avec un diagonale de $T$ contenant $*$. On va numéroter les
  diagonales de $T$ contenant $*$ par les éléments de $\extr$.

  Dans la définition de $\BB_T$ comme produit sur les diagonales de
  $T$, on peut séparer les contributions des diagonales strictement
  contenues dans les différents morceaux et la contribution des
  diagonales de $T$ séparant les morceaux. On va traiter séparément le
  morceau le plus à gauche et les autres morceaux. Ces autres morceaux
  sont numérotés par l'élément de $\extr$ qui les borde sur leur
  gauche.

  La contribution des diagonales entre les morceaux est 
  \begin{equation}
    \prod_{k \in \extr} \left( \x_{k+1}-\x_{n+1} \right).
  \end{equation}
  
  Considérons le premier morceau et soit $k_{\min}$ le plus petit
  élément de $\extr$. La contribution du premier morceau est
  \begin{equation}
    \prod_{{1\leq i\leq j< k_{\min}}\atop{(i,j) \in T}} 
    \left(\x_{i}-\x_{j+1}\right).
  \end{equation}

  Considérons maintenant $k\in \extr$ et le morceau correspondant,
  situé à droite de $k$. Soit $k'$ l'élément suivant de $\extr$ ou
  bien posons $k'=n+1$ si $k$ est le plus grand élément de $\extr$. La
  contribution du morceau $k$ est alors
  \begin{equation}
    \prod_{{k+1\leq i < k'}\atop{(k+1,i) \in T}} \left(\x_{k+1}-\x_{i+1}\right)
    \prod_{{k+1< i\leq j< k'}\atop{(i,j) \in T}} \left(\x_{i}-\x_{j+1}\right),
  \end{equation}
  où le premier facteur est associé aux diagonales du morceau $k$ qui
  contiennent le sommet $k$.

  On a donc montré que $\BB_T$ est le produit de facteurs associés à
  chaque morceau : pour le premier morceau,
  \begin{equation}
    \prod_{{1\leq i\leq j< k_{\min}}\atop{(i,j) \in T}} \left(\x_{i}-\x_{j+1}\right)
  \end{equation}
  et, pour le morceau à droite de $k$ dans $\extr$,
  \begin{equation}
    \left(\x_{k+1}-\x_{n+1} \right)
    \prod_{{k+1\leq i < k'}\atop{(k+1,i) \in T}} \left(\x_{k+1}-\x_{i+1}\right)
    \prod_{{k+1< i\leq j< k'}\atop{(i,j) \in T}} \left(\x_{i}-\x_{j+1}\right).
  \end{equation}

  Regardons maintenant l'image $D(T)$ de $T$ par la bijection. C'est
  la concaténation des images des morceaux de $T$. Par définition du
  monôme associé, celui-ci est le produit des contributions de chaque
  morceau avec un décalage des indices convenable et des contributions
  des pas verticaux initiaux des morceaux (sauf le premier).

  Par hypothèse de récurrence, la contribution du premier morceau est
  \begin{equation}
    \prod_{{1\leq i\leq j< k_{\min}}\atop{(i,j) \in T}} \x_{i}.
  \end{equation}

  La contribution du morceau entre $k\in \extr$ et l'élément suivant
  $k'$ de $\extr$ est donnée, par hypothèse de récurrence, par
  \begin{equation}
    \x_{k+1}^{\ell_k}    
    \prod_{{k+1< i\leq j< k'}\atop{(i,j) \in T}} \x_{i},
  \end{equation}
  où $\ell_k$ est le nombre de pas verticaux initiaux du morceau $k$.

  Par le lemme \ref{premierspas} appliqué au morceau $k$, on sait que
  le nombre $\ell_k$ de pas verticaux initiaux dans le morceau $k$ de
  $D(T)$ est égal à $1$ plus le nombre de diagonales dans le morceau
  $k$ de $T$ qui contiennent le sommet $k$. La contribution du morceau
  $k$ au monôme $M_{D(T)}$ est donc
  \begin{equation}
    \x_{k+1}\prod_{{k+1\leq i < k'}\atop{(k+1,i) \in T}} \x_{k+1}     
    \prod_{{k+1< i\leq j< k'}\atop{(i,j) \in T}} \x_{i}.
  \end{equation}
  
  On vérifie que le terme dominant de la contribution de chaque
  morceau à $\BB_T$ est bien égal à la contribution de chaque morceau
  à $M_{D(T)}$. En prenant le produit des contributions des morceaux,
  on obtient l'égalité voulue.
\end{proof}

\begin{theorem}
  \label{principal}
  Les polynômes $\BB_T$ associés aux triangulations forment une base
  de l'algèbre des coinvariants quasi-symétriques. Cette base est
  stable par le renversement des variables $\x_{i} \mapsto
  \x_{n+2-i}$. Les deux choix naturels d'ordre total sur les variables
  donnent deux bases monomiales, en prenant les monômes dominants des
  polynômes $\BB_T$.
\end{theorem}

\begin{proof}
  Dans \cite{bba}, il est démontré que les classes des monômes $M_D$
  associés aux chemins de Dyck forment une base de l'anneau des
  coinvariants quasi-symétriques. On déduit alors de la proposition
  \ref{correspond} que les classes des polynômes $\BB_T$ forment aussi
  une base. Le fait que cette base soit stable par le renversement est
  immédiat : l'image de $\BB_T$ est $\BB_{T'}$ où la triangulation
  $T'$ est obtenue par renversement de $T$. Enfin la dernière
  assertion est juste une reformulation de la proposition
  \ref{correspond} et son image par le renversement.
\end{proof}

\bibliographystyle{plain}
\bibliography{quasisym}

\end{document}